\input amstex
\magnification=\magstep1
\baselineskip=13pt
\documentstyle{amsppt}
\vsize=8.7truein
\NoRunningHeads
\def\aff{\operatorname{aff}}
\def\lin{\operatorname{lin}}
\def\vl{\operatorname{vol}}
\def\inte{\operatorname{int}}
\def\xx{\bold{x}}

\title Computing the Ehrhart quasi-polynomial of a rational simplex   \endtitle
\author Alexander Barvinok \endauthor
\address Department of Mathematics, University of Michigan, Ann Arbor,
MI 48109-1043, USA \endaddress
\email barvinok$\@$umich.edu  \endemail
\thanks This research was partially supported by NSF Grant DMS 0400617.
\endthanks
\abstract We present a polynomial time algorithm to compute any fixed 
number of the highest coefficients of the Ehrhart quasi-polynomial of a rational simplex.
Previously such algorithms were known for integer simplices and for rational 
polytopes of a fixed dimension. The algorithm is based on the formula relating the 
$k$th coefficient of the Ehrhart quasi-polynomial of a rational polytope to volumes of 
sections of the polytope by affine lattice subspaces parallel to $k$-dimensional faces 
of the polytope.  We discuss possible extensions and open questions. \endabstract
\date April 2005 \enddate
\endtopmatter
\document
 
 \head 1. Introduction and main results \endhead
 
 Let $P \subset {\Bbb R}^d$ be a rational polytope, that is, a polytope with rational 
 vertices. Let $t \in {\Bbb N}$ be a positive 
 integer such that the dilated polytope
 $$t P=\Bigl\{tx: \quad x \in P\Bigr\}$$ has integer vertices. As is known, see, for example,
 Section 4.6 of \cite{St97}, there 
 exist functions $e_i(P; \cdot ): {\Bbb N} \longrightarrow {\Bbb Q}$, $i=0, \ldots, d$, such that 
 $$e_i(P; n+t)=e_i(P; n) \quad \text{for all} \quad n \in {\Bbb N}$$
and
$$|nP \cap {\Bbb Z}^d| =\sum_{i=0}^d e_i(P; n) n^i \quad \text{for all} \quad n \in {\Bbb N}.$$
The function on the right hand side is called the {\it Ehrhart quasi-polynomial} of $P$.
It is clear that if $\dim P=d$ then $e_d(P; n)=\vl P$. In this paper, we are interested in 
 the computational complexity of the coefficients $e_i(P; n)$. 
 
 If the dimension $d$ is fixed in advance, the values of $e_i(P; n)$ for any 
 given $P$, $n$, and $i$ can be computed in 
 polynomial time by interpolation, as implied by a polynomial time algorithm to count 
 integer points in a polyhedron of a fixed dimension \cite{B94a}, \cite{BP99}.
  
 If the dimension $d$ is allowed to vary, it is an NP-hard problem to check whether 
 $P \cap {\Bbb Z}^d \ne \emptyset$, let alone to count integer points in $P$. This is true 
 even when $P$ is a rational simplex, as exemplified by the knapsack problem see, for 
 example, Section 16.6 of \cite{Sc86}.
 If the polytope $P$ is integral then the coefficients $e_i(P; n)=e_i(P)$ 
 do not depend on $n$. In that case, for any $k$ fixed in advance, computation of the 
 Ehrhart coefficient $e_{d-k}(P)$ reduces in polynomial time to computation of the 
 volumes of  the $(d-k)$-dimensional faces of $P$ \cite{B94b}. The algorithm is based
 on efficient formulas relating $e_{d-k}(P)$, volumes of the $(d-k)$-dimensional faces, 
 and cones of feasible directions at those faces, see \cite{Mo93}, \cite{BP99},
 and \cite{PT04}. In particular, if $P=\Delta$ is an integer simplex, there 
 is a polynomial time algorithm for computing $e_{d-k}(\Delta)$ as long as $k$ fixed in advance.
 
 In this paper, we extend the last result to {\it rational} simplices.
\bigskip
$\bullet$ Let us fix an integer $k \geq 0$. The paper presents a polynomial time algorithm, 
which, given an integer $d \geq k$, a rational simplex $\Delta \subset {\Bbb R}^d$, and a positive integer $n$, computes the value of $e_{d-k}(\Delta; n)$.
 \bigskip
 We present the algorithm in Section 7 and discuss its possible extensions  in Section 8.

In contrast to the case of an integral polytope, for a general rational polytope $P$ computation of
 $e_i(P; n)$ cannot be reduced 
 to computation of the volumes of faces and some functionals of the ``angles'' (cones of feasible direction) at the faces. A general result of McMullen \cite{Mc78}, see also \cite{MS83} and \cite{Mc93},
  asserts that the 
 contribution of the $i$-dimensional face $F$ of a rational polytope $P$ to the coefficient 
 $e_i(P; n)$ is a function of the volume of $F$, the cone of feasible directions of $P$ at $F$, {\it and} 
 the translation class of the affine hull $\aff(F)$ of $F$ modulo ${\Bbb Z}^d$. 
 \bigskip
 Our algorithm is based on a new structural result, Theorem 1.3 below, relating 
 the coefficient $e_{d-k}(P;n)$ to volumes of sections of $P$ by 
 affine lattice subspaces  parallel to faces $F$ of $P$ with $\dim F \geq d-k$.  Theorem 1.3 may 
 be of interest in its own right.
 
 \subhead (1.1) Valuations and polytopes \endsubhead 
 Let $V$ be a $d$-dimensional real vector space and let $\Lambda \subset V$ be a lattice, 
 that is, a discrete additive subgroup which spans $V$. A polytope $P \subset V$ is called
 a $\Lambda$-{\it polytope} or a {\it lattice polytope} if the vertices of $P$ belong to $\Lambda$.
 A polytope $P \subset V$ is called $\Lambda$-{\it rational} or just {\it rational} if 
 $tP$ is a lattice polytope for some positive integer $t$.
 
 For a set $A \subset V$, let $[A]: V \longrightarrow {\Bbb R}$ be the indicator of $A$:
 $$[A](x)=\cases 1 &\text{if\ } x \in A \\ 0 & \text{if\ } x \notin A. \endcases$$
 A complex-valued function $\nu$ on rational polytopes $P \subset V$ is called a {\it valuation} if 
 it preserves linear relations among indicators of rational polytopes:
 $$\sum_{i \in I} \alpha_i [P_i] =0 \Longrightarrow \sum_{i \in I} \alpha_i \nu(P_i)=0,$$
 where $P_i \subset V$ is a finite family of rational polytopes and $\alpha_i$ are rational 
 numbers. We consider only $\Lambda$-{\it valuations} or {\it lattice valuations} $\nu$ that satisfy
 $$\nu(P+u)=\nu(P) \quad \text{for all} \quad u \in \Lambda,$$
 see \cite{MS83} and \cite{Mc93}.

 A general result of McMullen \cite{Mc78} states that if 
$\nu$ is a lattice valuation, 
$P \subset V$ is a rational polytope, and $t \in {\Bbb N}$ is a number 
such that $tP$ is a lattice polytope then there exist functions 
$\nu_i(P; \cdot): {\Bbb N} \longrightarrow {\Bbb C}$, $i=0, \ldots, d$, such that 
$$\nu(nP)=\sum_{i=0}^d \nu_i(P; n) n^i \quad \text{for all} \quad n \in {\Bbb N}$$ 
and 
$$\nu_i(P; n+t)=\nu_i(P;  n) \quad \text{for all} \quad  n \in {\Bbb N}.$$
Clearly, if we compute $\nu(mP)$ for $m=n, n+t, n+2t, \ldots, n+dt$, we 
can obtain $\nu_i(P; n)$ by interpolation.

We are interested in the counting valuation $E$, where $V={\Bbb R}^d$, $\Lambda={\Bbb Z}^d$,  and
$$E(P)=|P \cap {\Bbb Z}^d|$$
is the number of lattice points in $P$.

The idea of the algorithm is to replace valuation $E$ by some other valuation, so that the 
coefficients $e_d(P; n), \ldots, e_{d-k}(P; n)$ remain intact, but the new valuation can be 
computed in polynomial  time on any given rational simplex $\Delta$, so that the desired coefficient
$e_{d-k}(\Delta; n)$ can be obtained by interpolation. 

 \subhead (1.2) Valuations $E_L$ \endsubhead
 Let $L \subset {\Bbb R}^d$ be a lattice subspace, that is, a subspace spanned by the points 
 $L \cap {\Bbb Z}^d$. Suppose that $\dim L=k$ and let $pr: {\Bbb R}^d \longrightarrow L$ be the orthogonal projection onto $L$. Let $P \subset {\Bbb R}^d$ be a rational polytope, let 
 $Q=pr(P)$, $Q \subset L$, be its projection, and let $\Lambda=pr\left({\Bbb Z}^d\right)$.
 Since $L$ is a lattice subspace, $\Lambda \subset L$ is a lattice.
 
  Let $L^{\bot}$ be the orthogonal complement of $L$. Then $L^{\bot} \subset {\Bbb R}^d$ is 
  a lattice subspace. We introduce the volume form $\vl_{d-k}$ on $L^{\bot}$ 
  which differs from the volume form inherited from ${\Bbb R}^d$ by a scaling factor 
  chosen so that the determinant of the lattice ${\Bbb Z}^d \cap L^{\bot}$ is 1.
   Consequently, the same volume form $\vl_{d-k}$ is carried by all 
  translations $x+L^{\bot}$, $x \in {\Bbb R}^d$.
  
 We consider the following quantity
 $$E_L(P)=\sum_{m \in \Lambda} \vl_{d-k}\left(P \cap \left(m+L^{\bot}\right)\right)=
 \sum_{m \in Q \cap \Lambda}  \vl_{d-k}\left(P \cap \left(m+L^{\bot}\right)\right)$$
 (clearly, for $m \notin Q$ the corresponding terms are 0).
 
 In words: we take all lattice translates of $L^{\bot}$, select those 
 that intersect $P$ and add the volumes of the intersections.
   
 Clearly, $E_L$ is a lattice valuation, so 
 $$E_L(nP)=\sum_{i=0}^d e_i(P, L; n) n^i$$
 for some periodic functions $e_i(P, L; \cdot)$. If $tP$ is an integer polytope for some $t \in {\Bbb N}$
 then 
 $$e_i(P, L; n+t)=e_i(P, L; n) \quad \text{for all} \quad n \in {\Bbb N}$$
 and $i=0, \ldots, d$.
 
 Note that if $L=\{0\}$ then $E_L(P)=\vl P$ and if $L={\Bbb R}^d$ then 
 $E_L(P)=| P \cap {\Bbb Z}^d|$, so the valuations $E_L$ interpolate between the volume 
 and the number of lattice points as $\dim L$ grows. 
 
 We prove that $e_{d-k}(P; n)$ can be represented as a linear combination of 
 
 \noindent $e_{d-k}(P, L; n)$ for some lattice subspaces $L$ with $\dim L \leq k$.
  \proclaim{(1.3) Theorem} Let us fix an integer $k\geq 0$. 
 Let $P \subset {\Bbb R}^d$ be a full-dimensional rational polytope
 and let $t$ be a positive integer such that $tP$ is an integer polytope. 
 For a $(d-k)$-dimensional face $F$ of $P$ let $\lin(F) \subset {\Bbb R}^d$ be the 
 $(d-k)$-dimensional subspace parallel to the affine hull $\aff(F)$ of $F$ and let 
  $L^F =\left(\lin F\right)^{\bot}$ be its orthogonal complement, so $L^F \subset {\Bbb R}^d$ 
  is a $k$-dimensional lattice subspace.
  
  Let ${\Cal L}$ be a finite collection of lattice subspaces which contains 
  the subspaces $L^F$ for  
  all $(d-k)$-dimensional faces $F$ of $P$ and is closed under intersections.
  For $L \in {\Cal L}$ let $\mu(L)$ be integer numbers such that the identity
  $$\left[ \bigcup_{L \in {\Cal L}}  L\right] =\sum_{L \in {\Cal L}} \mu(L) [L]$$
  holds for the indicator functions of the subspaces from ${\Cal L}$.
 
 Let us define
 $$\nu(nP)=\sum_{L \in {\Cal L}} \mu(L) E_{L}(nP) \quad \text{for} \quad n \in {\Bbb N}.$$
 Then there exist functions $\nu_i(P; \cdot): {\Bbb N} \longrightarrow {\Bbb Q}$, $i=0, \ldots, d$, such that  
 $$\nu(nP)=\sum_{i=0}^d \nu_i(P; n) n^i \quad \text{for all} \quad n \in {\Bbb N}, \tag1.3.1$$
 $$\nu_i(P; n+t)=\nu_i(P; n) \quad \text{for all} \quad n \in {\Bbb N}, \tag1.3.2$$
 and  
  $$e_{d-i}(P; n)=\nu_{d-i}(P; n) \quad \text{for all} \quad n \in {\Bbb N} \quad \text{and} \quad
  i=0, \ldots, k. \tag1.3.3$$
 \endproclaim
 
 We prove Theorem 1.3 in Section 4 after some preparations in Sections 2 and 3.
 \bigskip
 The advantage of working with valuations $E_L$ is that they are more amenable 
 to computations.
 \bigskip
 $\bullet$ Let us fix an integer $k \geq 0$. We present a polynomial time algorithm, which, given 
 an integer $d \geq k$, a $d$-dimensional rational simplex $\Delta \subset {\Bbb R}^d$,
 and a lattice subspace $L \subset {\Bbb R}^d$ such that $\dim L \leq k$, computes 
 $E_L(\Delta)$.
  \bigskip
 We present the algorithm in Section 6 after some preparations in Section 5.
 
 \subhead (1.4) The main ingredient of the algorithm to compute $e_{d-k}(\Delta; n)$ \endsubhead
 
 Theorem 1.3 allows us to reduce computation of $e_{d-k}(\Delta; n)$ to that 
 of $E_L(\Delta)$, where $L \subset {\Bbb R}^d$ is a lattice subspace and $\dim L \leq k$.
 Let us choose a particular lattice subspace $L$ with $\dim L=j \leq k$.
 
If $P=\Delta$ is a simplex, then the description of the orthogonal projection 
$Q=pr(\Delta)$ onto $L$ can be 
computed in polynomial time. Moreover, one can compute in polynomial time 
a decomposition of $Q$ into a union of non-intersecting polyhedral pieces $Q_i$, such that
$\vl_{d-j}\left( pr^{-1}(x)\right)$ is a polynomial on each piece $Q_i$. Thus computing of $E_L(\Delta)$ 
reduces to computing of the sum
$$\sum_{m \in Q_i \cap \Lambda} \phi(m),$$
where $\phi$ is a polynomial with $\deg \phi =d-j$, $Q_i \subset L$ is a polytope with 
$\dim Q_i =j \leq k$ and $\Lambda \subset L$ is a lattice. The sum is computed by applying the technique 
of ``short rational functions'' for lattice points in polytopes of a fixed dimension, cf.
\cite{BW03}, \cite{BP99}, and \cite{D+04}. 
 
 The algorithm for computing the sum of a polynomial over integer points in a polytope 
 is discussed in Section 5.
  
 \head 2. The Fourier expansions of $E$ and $E_L$ \endhead

Let $V$ be a $d$-dimensional real vector space with the scalar product $\langle \cdot, \cdot \rangle$
and the corresponding Euclidean norm $\| \cdot\|$. Let $\Lambda \subset V$ be a lattice
and let  $\Lambda^{\ast} \subset V$ be the {\it dual} or the {\it reciprocal} lattice
$$\Lambda^{\ast}=\Bigl\{x \in V: \quad \langle x, y \rangle \in {\Bbb Z} \quad 
\text{for all} \quad y \in \Lambda \Bigr\}.$$ 
 
For $\tau >0$, we introduce the {\it theta function}
$$\split \theta_{\Lambda} (x, \tau)=&\tau^{d/2} \sum_{m \in \Lambda} \exp\left\{ - \pi \tau \|x -m \|^2\right\}\\=
&(\det \Lambda)^{-1} \sum_{l \in \Lambda^{\ast}} \exp\left\{ - \pi \| l \|^2/\tau + 2\pi i \langle l, x \rangle 
\right\}, \quad \text{where} \quad x \in V. \endsplit$$
The last inequality is the reciprocity relation for theta series (essentially, the Poisson summation 
formula), see, for example, Section 69 of \cite{Be61}.

For a polytope $P$, let $\inte P$ denote the relative interior of $P$ and let 
$\partial P=P \setminus \inte P$ be the boundary of $P$.
 
\proclaim{(2.1) Lemma} Let $P \subset V$ be a full-dimensional polytope such that
$\partial P \cap \Lambda = \emptyset$. 
Then 
$$\split |P \cap \Lambda|=&\lim_{\tau \longrightarrow +\infty} \int_P \theta_{\Lambda}(x, \tau) \ dx \\= 
&(\det \Lambda)^{-1} \lim_{\tau \longrightarrow +\infty} \sum_{l \in \Lambda^{\ast}}
\exp\left\{ -\pi \| l\|^2/\tau \right\} \int_P \exp\{ 2 \pi i \langle l, x \rangle \} \ dx. \endsplit $$ 
\endproclaim
\demo{Proof}
As is known (cf., for example, Section B.5 of \cite{La02}), as $\tau \longrightarrow +\infty$,
the function $\theta_{\Lambda}(x, \tau)$ converges in the sense of 
distributions to the sum of the delta-functions concentrated at the points $m \in \Lambda$.
Therefore, for every smooth function $\phi: {\Bbb R}^d \longrightarrow {\Bbb R}$ with a compact
support, we have
$$\lim_{\tau \longrightarrow +\infty} \int_{{\Bbb R}^d} \phi(x) \theta_{\Lambda}(x, \tau) \ dx=
\sum_{m \in \Lambda} \phi(m).\tag2.1.1$$

Since $\partial P \cap {\Lambda}=\emptyset$, we can replace $\phi$ by the indicator function 
$[P]$ in (2.1.1).
{\hfill \hfill \hfill} \qed
\enddemo
\remark{Remark} If $\partial P \cap \Lambda \ne \emptyset$, the limit still exists but then it counts every lattice point $m \in \partial P$ with the weight equal to the ``solid angle'' of $m$ at $P$,
since every term $\exp\left\{-\pi \tau \|x-m\|^2\right\}$ is spherically symmetric about $m$.
 This connection between the solid angle valuation and the theta function was described 
by the author in the unpublished paper \cite{Ba92} (the paper is very different from 
paper \cite{B94b} which has the same title) and independently discovered 
by Diaz and Robins \cite{DR94}. Diaz and Robins used a similar approach based 
on Fourier analysis to express coefficients of the Ehrhart polynomial of an integer 
polytope in terms of cotangent sums \cite{DR97}. Banaszczyk \cite{B93a}  
obtained asymptotically optimal bounds in transference theorems for lattices by using 
a similar approach with theta functions, with the polytope $P$ replaced by a Euclidean
ball.

The formula of Lemma 2.1 can be considered as the Fourier expansion of the counting valuation.
\endremark

We need a similar result for valuation $E_L$ defined in Section 1.2. 

\proclaim{(2.2) Lemma} Let $P \subset {\Bbb R}^d$ be a 
full-dimensional polytope and let $L \subset {\Bbb R}^d$ be a 
lattice subspace with $\dim L=k$. Let $pr: {\Bbb R}^d \longrightarrow L$
be the orthogonal projection onto $L$, let $Q=pr(P)$, and let $\Lambda=pr({\Bbb Z}^d)$,
so $\Lambda \subset L$ is a lattice in $L$. 
Suppose that $\partial Q \cap  \Lambda =\emptyset$.

Then  
 $$E_L(P)=\lim_{\tau \longrightarrow +\infty} \sum_{l \in L \cap {\Bbb Z}^d } 
 \exp\left\{-\pi \|l\|^2/\tau\right\} \int_P \exp\{2 \pi i \langle l, x \rangle \} \ dx.$$
\endproclaim  

\demo{Proof} We observe that $L \cap {\Bbb Z}^d =\Lambda^{\ast}$.
For a vector $x \in {\Bbb R}^d$, let $x_L$ be the orthogonal projection of $x$ onto $L$.
Applying the reciprocity relation for theta functions in $L$, we 
write
 $$\split &\sum_{l \in L \cap {\Bbb Z}^d } 
  \exp \left\{ -\pi \|l\|^2/\tau + 2\pi i \langle l, x \rangle \right\} \\=
 &\sum_{l \in L \cap {\Bbb Z}^d }  \exp \left\{ -\pi \| l \|^2/\tau + 2\pi i \langle l, x_L \rangle \right\}\\=
 &(\det \Lambda) \tau^{k/2} \sum_{m \in \Lambda} \exp\left\{ - \pi \tau \|x_L - m\|^2 \right\}.
 \endsplit$$
As is known (cf., for example, Section B.5 of \cite{La02}), as $\tau \longrightarrow +\infty$,
the function 
$$g_{\tau}(x) =\tau^{k/2} \sum_{m \in \Lambda} \exp\left\{-\pi \tau \|x_L -m \|^2 \right\}$$
converges in the sense of distributions to the sum of the delta-functions concentrated 
on the subspaces $m +L^{\bot}$  (this is the set of points where $x_L=m$) for $m \in \Lambda$.

Therefore, for every smooth function $\phi: {\Bbb R}^d \longrightarrow {\Bbb R}$ with a
compact support, we have
$$\lim_{\tau \longrightarrow +\infty} \int_{{\Bbb R}^d} \phi(x) g_{\tau}(x) \ dx =
\sum_{m \in \Lambda} \int_{m+L^{\bot}} \phi(x) \ d_{L^{\bot}} x, \tag2.2.1$$
where $d_{L^{\bot}} x$ is the Lebesgue measure on $m+L^{\bot}$ induced from ${\Bbb R}^d$.  

Since $\partial Q \cap \Lambda =\emptyset$,  
each subspace $m + L^{\bot}$ for $m \in \Lambda$ 
either intersects the interior of $P$ or is at least some distance $\epsilon=\epsilon(P, L)>0$ 
away from $P$. Hence we may replace $\phi$ by the indicator $[P]$ in (2.2.1).

 Recall from Section 1.2 that measuring volumes in $m+L^{\bot}$, we scale the 
 volume form in $L^{\bot}$ induced from ${\Bbb R}^d$ so that the determinant of the lattice 
 $L^{\bot} \cap {\Bbb Z}^d$ is 1.
 One can observe that $\det \Lambda$ provides the required normalization factor, so
 $$(\det \Lambda) \int_{m+L^{\bot}} [P](x) \ d_{L^{\bot}}(x)=
 \vl_{d-k}\left( P \cap \left(m+L^{\bot}\right)\right).$$ 
 The proof now follows.
{\hfill \hfill \hfill} \qed
\enddemo

\remark{Remark} If $\partial Q \cap \Lambda \ne \emptyset$ the limit still 
exists but then for $m \in \partial Q \cap \Lambda$ the volume 
$\vl_{d-k}\left(P \cap \left(m+L^{\bot}\right)\right)$ is counted with the weight  
defined as follows: we find the minimal (under inclusion) face $F$ of $P$ such that 
$m+L^{\bot}$ is contained in $\aff(F)$ and let the weight equal to the solid angle of 
$P$ at $F$.
\endremark

 \head 3. Exponential valuations \endhead
 
Let $V$ be a $d$-dimensional Euclidean space, let $\Lambda \subset V$ be a lattice
and let $\Lambda^{\ast}$ be the reciprocal lattice.
 Let us choose a vector $l \in \Lambda^{\ast}$ and let us consider the integral
 $$\Phi_l(P)=\int_P \exp\{ 2 \pi i \langle l, x \rangle \} \ dx,$$
 where $dx$ is the Lebesgue measure in $V$.
 Note that for $l=0$ we have $\Phi_l(P)=\Phi_0(P)=\vl P$. 
 We have
 $$\Phi_l(P+a)=\exp\left\{ 2 \pi i \langle l, a \rangle \right\} \Phi_l(P) \quad \text{for all} \quad a \in V.$$
 It follows that  $\Phi_l$ is a $\Lambda$-valuation on rational polytopes $P \subset V$.
 
 If $l \ne 0$ then the following lemma (essentially, Stokes' formula) shows that 
 $\Phi_l$ can be expressed as a linear combination of exponential valuations on the facets of $P$.
 The proof can be found, for example, in \cite{B93b}.  
 
 \proclaim{(3.1) Lemma} Let $P \subset V$ be a full-dimensional polytope.
 For a facet $\Gamma$ of $P$, let $d_{\Gamma} x$ be the Lebesgue measure on 
 $\aff(\Gamma)$, and let $p_{\Gamma}$ be the unit outer normal to $\Gamma$. Then, for 
 every $l \in V\setminus 0$, we have
 $$\int_P \exp \left\{ 2 \pi i \langle l, x \rangle \right\} \ dx =\sum_{\Gamma} 
 {\langle l, p_{\Gamma} \rangle \over 2 \pi i \| l \|^2} 
\int_{\Gamma} \exp\left\{ 2 \pi i \langle l, x \rangle \right\}
 \ d_{\Gamma} x,$$
 where the sum is taken over all facets $\Gamma$ of $P$.
 \endproclaim
 
 Let $F \subset P$ be an $i$-dimensional face of $P$. Recall that by $\lin(F)$ we denote the 
 $i$-dimensional subspace of ${\Bbb R}^d$ that is parallel to the affine hull $\aff(F)$ of $F$.
 We need the following result. 
 
 \proclaim{(3.2) Theorem} Let $P \subset V$ be a rational full-dimensional polytope and let $t $ 
 be a positive integer such that $tP$ is a lattice polytope. Let 
 $\epsilon \geq  0$ be a rational number and let $a \in V$ be a vector. Let us choose $l \in \Lambda^{\ast}$. Then 
 there exist functions $f_i(P, \epsilon, a, l; \cdot): {\Bbb N} \longrightarrow {\Bbb C}$,
 $i=0, \ldots, d$, such that 
  $$\Phi_l\bigl((n+\epsilon)P+a \bigr)=\sum_{i=0}^d f_i(P, \epsilon, a, l; n) n^i \quad \text{for all} \quad n \in {\Bbb N} 
  \tag3.2.1$$
  and
  $$f_i(P, \epsilon,a, l; n+t) =f_i(P, \epsilon, a, l; n) \quad \text{for all} \quad n \in {\Bbb N}\tag3.2.2$$
 and $i=0,\ldots, d$.
 
 Suppose that $f_{d-k}(P, \epsilon, a, l; n) \ne 0$ for some $n$. Then there exists 
 a $(d-k)$-dimensional face $F$ of $P$ such that $l$ is orthogonal to $\lin(F)$.
 \endproclaim
  \demo{Proof} Since   
  $$\Phi_l(P+a)=\exp\left\{2 \pi i \langle l, a \rangle \right\} \Phi_l(P),$$
  without loss of generality we assume that $a=0$. We will denote 
  $f_i(P, \epsilon, 0, l; n)$ just by $f_i(P, \epsilon, l; n)$. 
    
  We proceed by induction on $d$. For $d=0$ the statement of the theorem
  obviously holds. Suppose that $d \geq 1$. If $l=0$ then 
 $\Phi_l\bigl((n+\epsilon)P\bigr)=(n+\epsilon)^d \vl P$ and the statement holds as well.
 
 Suppose that $l \ne 0$. For a facet $\Gamma$ of $P$, let 
 $\Lambda_{\Gamma}=\Lambda \cap \lin(\Gamma)$ and let 
 $l_{\Gamma}$ be the orthogonal projection of 
 $l$ onto $\lin(\Gamma)$. Thus $\Lambda_{\Gamma}$ is a lattice in the $(d-1)$-dimensional 
 Euclidean space $\lin(\Gamma)$ and $l_{\Gamma} \in \Lambda_{\Gamma}^{\ast}$, so we can define 
 valuations $\Phi_{l_{\Gamma}}$ on $\lin(\Gamma)$. Since $tP$ is a lattice polytope, for 
 every facet $\Gamma$ there is a vector $u_{\Gamma} \in V$ such that 
 $$\lin(\Gamma)=\aff(t \Gamma)-t u_{\Gamma} \quad \text{and} \quad t u_{\Gamma} \in \Lambda.$$
 
 Let $\Gamma'=\Gamma- u_{\Gamma}$, so $\Gamma' \subset \lin(\Gamma)$ is 
 a $\Lambda_{\Gamma}$-rational $(d-1)$-dimensional polytope such that 
 $t \Gamma'$ is a $\Lambda_{\Gamma}$-polytope. We have 
 $$ (n+\epsilon) \Gamma = (n+\epsilon)\Gamma' + (n+\epsilon) u_{\Gamma}.$$  
Applying Lemma 3.1 to $(n+\epsilon)P$, we get 
$$\Phi_l\bigl((n+\epsilon)P\bigr)=\sum_{\Gamma}  \psi(\Gamma, l; n) \Phi_{l_{\Gamma}}
\bigl((n+\epsilon) \Gamma' \bigr),$$
where 
$$\psi(\Gamma, l; n)={\langle l, p_{\Gamma} \rangle \over 2 \pi i \|l\|^2} 
\exp\bigl\{ 2 \pi i (n+\epsilon) \langle l, u_{\Gamma} \rangle \bigr\} $$
and the sum is taken over all facets $\Gamma$ of $P$.
 
 Since $t u_{\Gamma} \in \Lambda$ and $l \in \Lambda^{\ast}$, we have 
 $$\psi(\Gamma, l; n+t) =\psi(\Gamma, l; n) \quad \text{for all} \quad n \in {\Bbb N}.$$ 
  Hence, applying the induction hypothesis, we may write
  $$f_i(P, \epsilon, l; n)=\sum_{\Gamma} \psi(\Gamma, l; n) f_i(\Gamma', \epsilon, l_{\Gamma}; n)
  \quad \text{for all} \quad n \in {\Bbb N}$$
  and $i=0, \ldots, d-1$ and $f_d(P, \epsilon, l; n) \equiv 0$. Hence (3.2.1)-(3.2.2) 
  follows by the induction hypothesis. 
  
  If $f_{d-k}(P, \epsilon, l; n) \ne 0$ then there is a facet $\Gamma$ of $P$ such that 
  $f_{d-k}(\Gamma', \epsilon, l_{\Gamma}; n) \ne 0$. By the induction hypothesis, there is 
  a face $F'$ of $\Gamma'$ such that $\dim F'=d-k$, and $l_{\Gamma}$ is orthogonal to $\lin(F')$.
  Then $F=F'+u_{\Gamma}$ is a $(d-k)$-dimensional face of $P$, $\lin(F')=\lin(F)$ and 
  $l$ is orthogonal to $\lin(F)$, which completes the proof.
  {\hfill \hfill \hfill} \qed
  \enddemo
  
  \head 4. Proof of Theorem 1.3 \endhead
  
  First, we discuss some ideas relevant to the proof.
  \subhead (4.1) Shifting a valuation by a polytope \endsubhead
  Let $V$ be a $d$-dimensional real vector space, let $\Lambda \subset V$ be a lattice, and 
  let $\nu$ be a $\Lambda$-valuation on rational polytopes. 
  Let us fix a rational polytope $R\subset V$. McMullen \cite{Mc78} observed that 
  the function $\mu$ defined by 
  $$\mu(P)=\nu(P+R)$$
  is a $\Lambda$-valuation on rational polytopes $P$.
 Here ``+'' stands for the Minkowski sum:
$$P+R=\Bigl\{x+y: \quad x \in P, y \in R \Bigr\}.$$
This result follows since the transformation $P \longmapsto P+R$ preserves linear 
dependencies among indicators of polyhedra, cf. \cite{MS83}.

Let $t$ be a positive integer such that $tP$ is a lattice polytope.
McMullen \cite{Mc78} deduced that there 
exist functions $\nu_i(P, R; \cdot): {\Bbb N} \longrightarrow {\Bbb C}$, $i=0, \ldots, d$, such that
$$\nu\bigl(nP+R\bigr)=\sum_{i=0}^d \nu_i(P, R; n) n^i \quad \text{for all} \quad n \in {\Bbb N}$$
and
$$\nu_i(P, R; n +t)=\nu_i(P, R; n) \quad \text{for all} \quad n \in {\Bbb N}.$$
  
 \subhead (4.2) Continuity properties of valuations $E$ and $E_L$ \endsubhead 
  Let $R \subset {\Bbb R}^d$ be a full-dimensional rational polytope containing the origin
  in its interior. Then for 
  every polytope $P \subset {\Bbb R}^d$ and every $\epsilon >0$ we have $P \subset \bigl(P + \epsilon R\bigr)$.
  We observe that
  $$\big| (P + \epsilon R) \cap {\Bbb Z}^d \big|=
  |P \cap {\Bbb Z}^d|,$$
  for all sufficiently small $\epsilon >0$. If $P$ is a rational polytope, the supporting affine hyperplanes
  of the facets of $nP$ for $n \in {\Bbb N}$ are split among finitely many translation classes 
  modulo ${\Bbb Z}^d$. Therefore, there exists $\delta=\delta(P, R) >0$ such that 
  $$\big| (nP + \epsilon R) \cap {\Bbb Z}^d \big|=
  |nP \cap {\Bbb Z}^d| \quad \text{for all} \quad 0< \epsilon < \delta \quad \text{and all} \quad 
  n \in {\Bbb N}.$$ 
  
  We also note that for every rational subspace $L \subset {\Bbb R}^d$, we have 
  $$\lim_{\epsilon \longrightarrow 0+} E_L\bigl(P + \epsilon R\bigr) =E_L(P).$$
  We will use the perturbation $P \longmapsto P + \epsilon R$ to push valuations $E$ and 
  $E_L$ into a sufficiently generic position, so that we can apply Lemmas 2.1--2.2 without having
   to deal with various boundary effects. 
  This is somewhat similar in spirit to the idea of  \cite{BS05}.
  
  \subhead (4.3) Linear identities for quasi-polynomials \endsubhead
  Let us fix positive integers $t$ and $d$. Suppose that we have a possibly infinite family of 
  quasi-polynomials $p_l: {\Bbb N} \longrightarrow {\Bbb C}$ of the type
  $$p_l(n)=\sum_{i=0}^d p_i(l; n)n^i \quad \text{for all} \quad n \in {\Bbb N},$$
  where functions $p_i(l; \cdot): {\Bbb N} \longrightarrow {\Bbb C}$, $i=0, \ldots, d$, satisfy
  $$p_i(l; n)=p_i(l; n+t) \quad \text{for all} \quad n\in {\Bbb N}.$$ 
  Suppose further that $p: {\Bbb N} \longrightarrow {\Bbb C}$ is yet another quasi-polynomial
  $$p(n)=\sum_{i=0}^d p_i(n) n^i \quad \text{where} \quad p_i(n+t)=p_i(n) \quad \text{for all}
  \quad n \in {\Bbb N}.$$
  Finally, suppose that $c_l(\cdot): {\Bbb R}_+ \longrightarrow {\Bbb C}$ is a family of 
  functions and that 
  $$p(n)=\lim_{\tau \longrightarrow +\infty} \sum_l c_l(\tau) p_l(n) \quad \text{for all}
  \quad  n \in {\Bbb N}$$
  and that the series converges absolutely for every $n \in {\Bbb N}$ and every $\tau >0$.
  
  Then we claim that for $i=0, \ldots, d$ we have
   $$p_i(n)=\lim_{\tau \longrightarrow +\infty} \sum_l c_l(\tau) p_i(l; n) \quad \text{for all} 
   \quad n \in {\Bbb N}$$
  and that the series converges absolutely for every $n \in {\Bbb N}$ and every $\tau >0$.

  This follows since $p_i(n)$, respectively $p_i(l; n)$, can be expressed as 
  linear combinations of $p(m)$, respectively $p_l(m)$, for 
  $m=n, n+t, \cdots, n+ td$ with the coefficients depending on $m, n, t$, and $d$ only.  
  \bigskip
  Now we are ready to prove Theorem 1.3.
  \bigskip
   \subhead (4.4) Proof of Theorem 1.3 \endsubhead Let us fix 
    a rational polytope $P \subset {\Bbb R}^d$ as defined in the statement of the theorem. 
    For $L \in {\Cal L}$ let  $P_L \subset L$ be the orthogonal projection of 
    $P$ onto $L$ and let $\Lambda_L \subset L$ be the orthogonal projection of 
    ${\Bbb Z}^d$ onto $L$. 
    
    Let $a \in \inte P$ be a rational vector and let 
    $$R=P-a.$$
    Hence $R$ is a rational polytope containing the origin in its interior. Let $R_L$ denote the orthogonal projection
    of $R$ onto $L$.
    
    Since $P$ is a rational polytope and ${\Cal L}$ is a 
    finite set of rational subspaces, there exists $\delta=\delta(P, R)>0$ such that for all 
    $0<\epsilon < \delta$ and all $n \in {\Bbb N}$, we have 
    $$\bigl(nP + \epsilon R\bigr) \cap {\Bbb Z}^d =n P \cap {\Bbb Z}^d \quad \text{and} \quad 
    \partial\bigl(nP +\epsilon R\bigr) \cap {\Bbb Z^d}=\emptyset \quad \text{for all} \quad n \in {\Bbb N}
    \tag1$$
    and for all $L \in {\Cal L}$, we have
    $$\bigl(n P_L + \epsilon R_L \bigr) \cap \Lambda_L=n P_L \cap \Lambda_L \quad \text{and} 
    \quad \partial\bigl(nP_L + \epsilon R_L \bigr) \cap \Lambda_L =\emptyset \quad \text{for all} \quad n \in {\Bbb N},\tag2$$ 
    cf. Section 4.2.
  Let us choose any rational $0< \epsilon < \delta$.
  
  Because of (1), we can write 
  $$\big| (nP+ \epsilon R) \cap {\Bbb Z}^d\big|=\sum_{i=0}^d e_i(P; n)n^i \quad \text{for all}\quad n \in 
  {\Bbb N} \tag3 $$
  and by Lemma 2.1 we get 
  $$\big|(nP +\epsilon R) \cap {\Bbb Z}^d \big|=\lim_{\tau \longrightarrow +\infty} 
  \sum_{l \in {\Bbb Z}^d} \exp\{-\pi \|l\|^2/\tau\} \Phi_l\bigl(nP+\epsilon R\bigr), \tag4$$
  where $\Phi_l$ are the exponential valuations of Section 3. 
  Since $\Phi_l$ is a ${\Bbb Z}^d$-valuation, by Section 4.1 there exist functions
  $f_i(P,  \epsilon, l; \cdot) :{\Bbb N} \longrightarrow {\Bbb C}$, $i=0, \ldots, d$, such that 
  $$\Phi_l(nP + \epsilon R)=\sum_{i=0}^d f_i(P, \epsilon, l; n)n^i \quad \text{for} \quad n \in {\Bbb N}
  \tag5$$
  and
  $$f_i(P, \epsilon, l; n+t)=f_i(P, \epsilon, l; n) \quad \text{for all} \quad n \in {\Bbb N}.\tag6$$
  
  Moreover, we can write
  
  $$nP+\epsilon R=nP + \epsilon(P-a)=(n+\epsilon) P - \epsilon a.$$
  Therefore, by Theorem 3.2, for $i \leq k$ we have 
  $f_{d-i}(P, \epsilon, l; n)=0$ unless 
  $l \in L^F$ for some face $F$ of $P$ with  $\dim F =d-k$.
  
  Therefore, combining (3)--(6) and Section 4.3, we obtain for all $0 \leq i \leq k$ and all $n \in {\Bbb N}$
  $$\split e_{d-i}(P; n) =&\lim_{\tau \longrightarrow +\infty} 
  \sum_{l \in {\Bbb Z}^d} \exp\left\{-\pi \|l\|^2/\tau\right\} f_{d-i}(P, \epsilon, l; n) \\=
  &\lim_{\tau \longrightarrow +\infty} 
  \sum_{l \in \bigcup_{L \in {\Cal L}} \left( L \cap {\Bbb Z}^d\right)}
   \exp\left\{-\pi \|l\|^2/\tau\right\} f_{d-i}(P, \epsilon, l; n), \endsplit$$
  since vectors $l \in {\Bbb Z}^d$ outside of subspaces $L \in {\Cal L}$ contribute 0 to the sum.
  Therefore, for $0 \leq i \leq k$ and all $n \in {\Bbb N}$
  $$e_{d-i}(P; n) =\lim_{\tau \longrightarrow +\infty} \sum_{L \in {\Cal L}} \mu(L)
  \sum_{l \in L \cap {\Bbb Z}^d} \exp\left\{-\pi\|l\|^2/\tau\right\} f_{d-i}(P, \epsilon, l; n)
  \tag7$$
  On the other hand, because of (2), by Lemma 2.2 we get for all $L \in {\Cal L}$ and all $n \in {\Bbb N}$
  $$E_L\bigl(nP+\epsilon R\bigr) =\lim_{\tau \longrightarrow +\infty} \sum_{l \in L \cap {\Bbb Z}^d}  
  \exp\left\{-\pi \|l\|^2/\tau\right\} \Phi_l\bigl(nP+\epsilon R\bigr). \tag8$$
  Since $E_L$ are ${\Bbb Z}^d$-valuations, by Section 4.1 there exist functions 
  
  \noindent $e_i(P, \epsilon, L; \cdot): {\Bbb N} \longrightarrow  {\Bbb Q}$, $i=0, \ldots, d$, such that 
  $$E_L\bigl(nP+\epsilon R\bigr)=\sum_{i=0}^d e_i(P, \epsilon, L; n) n^i 
   \quad \text{for all} 
  \quad n \in {\Bbb N} \tag9$$
  and 
  $$e_i(P, \epsilon, L; n+t)=e_i(P, \epsilon, L; n) \quad \text{for all} \quad n \in {\Bbb N}. \tag10$$
Combining (5)-(6) and (8)--(10), by Section 4.3 we conclude
$$e_{d-i}(P, \epsilon, L; n) =\lim_{\tau \longrightarrow +\infty} \sum_{l \in L \cap {\Bbb Z}^d}  
  \exp\left\{-\pi \|l\|^2/\tau\right\} f_{d-i}(P, \epsilon, l; n) \quad \text{for all} \quad n \in {\Bbb N}.$$
  Therefore, by (7), for $0 \leq i \leq k$ we have
  $$e_{d-i}(P; n)=\sum_{L \in {\Cal L}} \mu(L) e_{d-i}(P, \epsilon, L; n) \quad \text{for all} \quad n \in {\Bbb N}.
  \tag11$$
  Since $E_L$ is a ${\Bbb Z}^d$-valuation, there exist functions 
  $e_i(P, L; \cdot): {\Bbb N} \longrightarrow {\Bbb Q}$, $i=0, \ldots, d$, such that  
  $$E_{L}(nP)=\sum_{i=0}^d e_i(P, L; n) n^i \quad \text{for all} \quad n \in {\Bbb N} \tag12$$
  and 
  $$e_i(P, L; n+t) =e_i(P, L; n) \quad \text{for all} \quad n \in {\Bbb N}.$$
  Let us choose an $m \in {\Bbb N}$. Substituting $n=m, m+t, \ldots, m+dt$ in (12), we obtain 
  $e_i(P, L; m)$ as a linear combination of $E_L(nP)$ with coefficients depending 
  on $n$, $m$, $t$, and $d$ only. Similarly, substituting $n=m, m+t, \ldots, m+dt$ in (9), we obtain 
  $e_i(P, \epsilon, L; m)$ as the same linear combination of 
  $E_{L}\bigl(nP+\epsilon R\bigr)$. Since volumes are continuous functions, in view of (2) (see
  also Section 4.2), we get
  $$\lim_{\epsilon \longrightarrow 0+} E_L\bigl(nP+\epsilon R\bigr)=E_L(nP)
  \quad \text{for} \quad n=m, m+t, \ldots, m+dt.$$
  Therefore,
  $$\lim_{\epsilon \longrightarrow 0+} e_i(P, \epsilon, L; m)=e_i(P, L; m) \quad \text{for all} \quad
  m \in {\Bbb N}.$$
  Taking the limit as $\epsilon \longrightarrow +0$ in (11), we obtain for $0 \leq i \leq k$
  $$e_{d-i}(P; n)=\sum_{L \in {\Cal L}} \mu(L) e_{d-i}(P, L; n) \quad \text{for all} \quad n \in {\Bbb N}.$$
  To complete the proof, we note that
  $$\nu_{d-i}(P, L; n)=\sum_{L \in {\Cal L}} \mu(L) e_{d-i}(P, L; n).$$
  {\hfill \hfill \hfill} \qed
 
 \head 5. Summing up a polynomial over integer points in a rational polytope \endhead
 
 Let us fix a positive integer $k$ and let us consider the following situation.
 Let $Q \subset {\Bbb R}^k$ be a rational polytope, let $\inte Q$ be the relative interior of $Q$ and 
 let $f: {\Bbb R}^k \longrightarrow {\Bbb R}$ be a polynomial with rational 
 coefficients. We want to compute the 
 value 
 $$\sum_{m \in \inte Q \cap {\Bbb Z}^k} f(m). \tag5.1$$
 We claim that as soon as the dimension $k$ of the polytope $Q$ is fixed, there is a polynomial
 time algorithm to do that. We assume that the polytope $Q$ is a given by the list of its 
 vertices and the polynomial $f$ is given by the list its coefficients.
  
 For an integer point $m=(\mu_1, \ldots, \mu_k)$, let 
 $$\xx^m=x_1^{\mu_1} \cdots x_k^{\mu_k} \quad \text{for} \quad \xx=(x_1, \ldots, x_k)$$
 be the Laurent monomial in $k$ variables $\xx=(x_1, \ldots, x_k)$. We use the following 
 result \cite{BP99}.
 \subhead (5.2) The short rational function algorithm \endsubhead 
  Let us fix $k$. There is a polynomial time algorithm, which,
 given a rational polytope $Q \subset {\Bbb R}^k$ computes 
 the generating function (Laurent polynomial)
 $$S(Q; \xx)=\sum_{m \in \inte Q \cap {\Bbb Z}^k} \xx^m$$ in the form
 $$S(Q;  \xx)=\sum_{i \in I} \epsilon_i {\xx^{a_i} \over (1- \xx^{b_{i1}}) \cdots (1-\xx^{b_{ik}})},$$
 where $a_i \in {\Bbb Z}^k$, $b_{ij} \in {\Bbb Z}^k \setminus \{0\}$ and $\epsilon_i \in {\Bbb Q}$.
 In particular, the number $|I|$ of fractions is bounded by a polynomial in the input size of $Q$.
 \bigskip
 
 Our first step is computing the generating function 
$$S(Q, f; \xx)=\sum_{m \in Q \cap {\Bbb Z}^k} f(m) \xx^m.$$
Our approach is similar to that of \cite{D+04}.
\subhead (5.3) The algorithm for computing $S(Q, f; \xx)$ \endsubhead
 We observe that 
 $$S(Q, f; \xx) = f\left( x_1 {\partial \over \partial x_1}, \ldots,
 x_k {\partial \over \partial x_k} \right) S(Q; \xx).$$
 We compute $S(Q; \xx)$ as in Section 5.2. 
 
 Let $a=(\alpha_1, \ldots, \alpha_k)$ be an integer vector, 
 let $b_j=(\beta_{j1},  \ldots, \beta_{jk})$ be non-zero integer vectors for 
 $j=1, \ldots, k$ and let $\gamma_1, \ldots, \gamma_k$ be positive integers.
 Then
 $$\split &\left(x_i {\partial \over \partial x_i} \right) {\xx^a \over (1-\xx^{b_1})^{\gamma_1} \cdots 
 (1-\xx^{b_k})^{\gamma_k}} \\=
& \alpha_i {\xx^a \over (1-\xx^{b_1})^{\gamma_1} \cdots (1-\xx^{b_k})^{\gamma_k}} + 
\sum_{j=1}^k \gamma_j \beta_{ji} 
 {\xx^{a+b_j} \over (1-\xx^{b_j})^{\gamma_j+1}} 
\prod_{s \ne j} {1 \over (1-\xx^{b_s})^{\gamma_s}} . \endsplit$$
 Consecutively applying the above formula and collecting similar fractions, we compute
 $$f\left( x_1{\partial \over \partial x_1}, \cdots, x_k {\partial \over \partial x_k} \right) 
 {\xx^a \over (1-\xx^{b_1}) \cdots (1-\xx^{b_k})}$$ 
 as an expression of the type
 $$\sum_j \rho_j { \xx^{a_j} \over (1-\xx^{b_1})^{\gamma_{j1}} \cdots (1-\xx^{b_k})^{\gamma_{jk}}},
 \tag5.3.1$$
 where $\rho_j \in {\Bbb Q}$, $\gamma_{j1}, \ldots, \gamma_{jk}$ are non-negative
 integers satisfying $\gamma_{j1} + \ldots + \gamma_{jk} \leq k+\deg f$ and 
 $a_j$ are vectors of the type 
 $$a_j=a+ \mu_1 b_1 + \ldots + \mu_k b_k,$$
 where $\mu_i$ are non-negative integers and $\mu_1 + \ldots + \mu_k \leq \deg f$.  
 The number of terms in (5.3.1) is bounded by $(\deg f)^{O(k)}$, which shows that for 
 a $k$ fixed in advance, the algorithm runs in polynomial time.
 
Consequently, $S(Q, f; \xx)$ is computed in polynomial time.
\bigskip
Formally speaking, to compute the sum (5.1), we have to substitute 
$x_i=1$ into the formula for $S(Q, f; \xx)$. This, however, cannot be done in a
straightforward way since $\xx=(1, \ldots, 1)$ is a pole of every fraction in the expression for 
$S(Q, f; \xx)$. Nevertheless, the substitution can be done via efficient computation of 
the relevant residue of $S(Q, f; \xx)$ as described in \cite{B94a} and 
\cite{BW03}.

\subhead (5.4) The algorithm for computing the sum \endsubhead
The output of Algorithm 5.3 represents $S(Q, f; \xx)$ in the general form 
$$S(Q, f; \xx) =\sum_{i \in I} \epsilon_i {\xx^{a_i} \over (1-\xx^{b_{i1}})^{\gamma_{i1}} \cdots 
(1-\xx^{b_{ik}})^{\gamma_{ik}}},$$
where $\epsilon_i \in {\Bbb Q}$, $a_i \in {\Bbb Z}^k$, $b_{ij} \in {\Bbb Z}^k \setminus \{0\}$, 
and $\gamma_{ij} \in {\Bbb N}$ such that $\gamma_{i1} + \ldots + \gamma_{ik} \leq k+\deg f$
for all $i \in I$.
 
 Let us choose a vector $l \in {\Bbb Q}^k$, 
 $l=(\lambda_1, \ldots, \lambda_k)$ such that $\langle l, b_{ij} \rangle \ne 0$ for 
 all $i,j$ (such a vector can be computed in polynomial time, cf. \cite{B94a}).
 For a complex $\tau$, let 
 $$\xx(\tau)=\left(e^{\tau \lambda_1}, \ldots, e^{\tau \lambda_k} \right).$$
 We want to compute the limit 
 $$\lim_{\tau \longrightarrow 0} G(\tau) \quad \text{for} \quad G(\tau)=S\bigl(Q, f; \xx(\tau)\bigr).$$
 In other words, we want to compute the constant term of the Laurent expansion of $G(\tau)$ 
 around $\tau=0$.
 
 Let us consider a typical fraction
 $${\xx^a \over (1-\xx^{b_1})^{\gamma_1}  \cdots (1-\xx^{b_k})^{\gamma_k}}.$$
 Substituting $\xx(\tau)$, we get the expression 
 $${e^{\alpha \tau} \over (1-e^{\tau \beta_1})^{\gamma_1} \cdots (1-e^{\tau \beta_k})^{\gamma_k}}, \tag5.4.1$$
  where  $\alpha=\langle a, l \rangle$ and $\beta_i=\langle b_i, l \rangle$ for $i=1, \ldots, k$.
 The order of the pole at $\tau=0$ is $D=\gamma_1 + \ldots + \gamma_k \leq k+\deg f$.
 To compute the constant term of the Laurent expansion of (5.4.1) at $\tau=0$, we do the following.
 
 We compute the polynomial
 $$q(\tau)=\sum_{i=0}^D {\alpha^i \over i!} \tau^i $$
 that is the truncation at $\tau^D$ of the Taylor series expansion of $e^{\alpha \tau}$.
 For $i=1, \ldots, k$ we compute the polynomial $p_i(\tau)$ with $\deg p_i =D$ such that
 $${\tau \over 1- e^{\tau \beta_i}} = p_i(\tau)  + \quad \text{terms of higher order in\ } \tau $$
 at $\tau=0$. Consecutively multiplying polynomials $\mod \tau^{D+1}$ we compute a 
 polynomial $u(\tau)$ with $\deg u= D$ such that 
 $$q(\tau) p_1^{\gamma_1}(\tau) \cdots p_k^{\gamma_k}(\tau) \equiv u(\tau) \mod \tau^{D+1}.$$
  The coefficient of $\tau^D$ in $u(\tau)$ is the desired constant term of the Laurent expansion.
 
\head 6. Computing $E_L(\Delta)$ \endhead

Let us fix a positive integer $k$.
 Let $\Delta \subset {\Bbb R}^d$ be a rational simplex given by the list of its 
vertices and let $L \subset {\Bbb R}^d$ be a rational subspace given its basis 
and such that $\dim L=k$. In this section, we describe a polynomial time algorithm for 
computing the value of $E_L(\Delta)$ as defined in Section 1.2.
 
 Let $pr: {\Bbb R}^d \longrightarrow L$ be the orthogonal projection. 
 We compute the vertices of the polytope $Q=pr(\Delta)$ and a basis of the lattice 
 $\Lambda=pr({\Bbb Z}^d)$. 
 For basic lattice algorithms see \cite{Sc93} and \cite{G+93}.
 
 As is known, as $x \in \Delta$ varies, the function 
 $$\phi(x)=\vl_{d-k} \left(P_x\right) \quad \text{where}  \quad P_x=\left(\Delta \cap \left(x + L^{\bot}\right)\right)$$
 is a piece-wise polynomial on $Q$. Our first step consists of computing a decomposition
 $$Q=\bigcup_i C_i \tag6.1$$
 such that $C_i \subset Q$ are rational polytopes (chambers) with pairwise disjoint interiors and 
 polynomials $\phi_i: L \longrightarrow {\Bbb R}$ such that $\phi_i(x)=\phi(x)$ for $x \in C_i$. 
 
 We observe that every vertex of $P_x$ is the intersection of 
 $x+L^{\bot}$ and some $k$-dimensional face $F$ of $\Delta$.
 
 For every face $G$ of $\Delta$ with $\dim G=k-1$ and such that $\aff(G)$ is not 
 parallel to $L^{\bot}$, let us compute 
 $$A_G=\Bigl\{x \in L: \quad x + L^{\bot} \cap \aff(G) \ne \emptyset \Bigr\}.$$
 Then $A_G$ is an affine hyperplane in $L$. Then number of different hyperplanes $A_G$ is
 $d^{O(k)}$ and hence they cut $Q$ into at most $d^{O(k^2)}$ polyhedral chambers 
 $C_i$, cf. Section 6.1 of \cite{Ma02}. As long as $x$ stays within the relative interior of a chamber $C_i$, the strong combinatorial type of $P_x$ does not change (the facets of $P_x$ move parallel to themselves) and hence the restriction $\phi_i$ of $\phi$ 
 onto $C_i$ is a polynomial, cf. Section 5.1 of  \cite{Sc93}. Since in the $(d-k)$-dimensional space 
 $x+L^{\bot}$ the polytope $P_x$ is defined by $d$ linear inequalities, $\phi_i$ can be 
 computed in polynomial time, see \cite{GK94} and \cite{B93b}.
 
 The decomposition 6.1 gives rise to the formula
 $$[Q]=\sum_j [Q_j],$$
where $Q_j$ are open faces of the chambers $C_i$ (the number of such faces is bounded 
by a polynomial in $d$), cf. Section 6.1 of \cite{Ma02}. 
Hence we have
$$E_L(\Delta)=\sum_j \sum_{m \in Q_j \cap \Lambda} \phi(m).$$
We compute inner sums as described in Section 5.  
 
 \head 7. Computing $e_{d-k}(\Delta; n)$ \endhead 
 
 Let us fix a an integer $k \geq 0$. We describe our algorithm, which, given 
 a positive integer $d \geq k$, a rational simplex $\Delta \subset {\Bbb R}^d$ (defined, for example, 
 by the list of its vertices), and a positive integer $n$, computes the number $e_{d-k}(\Delta; n)$.
 
We use Theorem 1.3.
\subhead (7.1) Computing the set ${\Cal L}$ of subspaces  \endsubhead 
 We compute subspaces $L$ and numbers $\mu(L)$ described in Theorem 1.3.
 Namely, for each $(d-k)$-dimensional face $F$ of $\Delta$, we compute a basis of the 
 subspace $L^F =(\lin F)^{\bot}$. Hence $\dim L^{F} \leq k$. Clearly, the number of 
 distinct subspaces $L^F$ is $d^{O(k)}$. We let ${\Cal L}$ be the set consisting 
 of the subspaces $L^F$ and all other subspaces obtained
 as intersections of $L^F$. We compute ${\Cal L}$ in $k$ (or fewer) steps. Initially, we 
 let 
 $${\Cal L}:=\Bigl\{L^F:\quad F \quad \text{is a} \quad (d-k)\text{-dimensional face of} \quad 
 \Delta \Bigr\}.$$
 Then, on every step, we consider the previously constructed 
 subspaces $L \in {\Cal L}$, consider the pairwise intersections $L\cap L^F$
 as $F$ ranges over the $(d-k)$-dimensional faces of $\Delta$ and add the obtained
 subspace $L \cap L^F$ to the set
 ${\Cal L}$ if it was not already there. If no new subspaces are obtained, we stop. Clearly, in the 
 end of this process, we will obtain all subspaces $L$ that are intersections of 
 different $L^{F_i}$. Since $\dim L^F = k$, each subspace 
 $L \in {\Cal L}$ is an intersection of some $k$ subspaces $L^{F_i}$. Hence the process stops 
 after $k$ steps and the total number $|{\Cal L}|$ of subspaces is $d^{O(k^2)}$.
 
 Having computed the subspaces $L \in {\Cal L}$, we compute the numbers $\mu(L)$ as follows.

For each pair of subspaces $L_1, L_2 \in {\Cal L}$ such that
$L_1 \subseteq L_2$, we compute the number $\mu(L_1, L_2)$
recursively: if $L_1=L_2$ we let $\mu(L_1, L_2)=1$. Otherwise, we let 
$$\mu(L_1, L_2)=-\sum \Sb L \in {\Cal L} \\ L_1 \subseteq L \subsetneq L_2 \endSb \mu(L_1, L).$$
 In the end, for each $L \in {\Cal L}$, we let
 $$\mu(L)=\sum \Sb L_1 \in {\Cal L} \\ L  \subseteq L_1  \endSb \mu(L, L_1).$$ 
Hence $\mu(L_i, L_j)$ are the values of the M\"obius function on 
 the set ${\Cal L}$ partially ordered by inclusion, so 
 $$\left[\bigcup_{L \in {\Cal L}} L\right] =\sum_{L \in {\Cal L}} \mu(L) [L]$$
 follows from the M\"obius inversion formula, cf. Section 3.7 of \cite{St97}.
 \bigskip
 Now, for each $L \in {\Cal L}$ and
 $m=n, n+t, \ldots , n+ td$ we compute the values of $E_L(m\Delta)$ as in Section 6, compute
 $$\nu(m\Delta)=\sum_{L \in {\Cal L}} \mu(L) E_L(m \Delta)$$
 and find $\nu_{d-k}(\Delta; n)=e_{d-k}(\Delta, n)$ by interpolation.
 
 \head 8. Possible extensions and further questions  \endhead
 
\subhead (8.1) Computing more general expressions \endsubhead 
Let $P \subset {\Bbb R}^d$ be a rational polytope, let $\alpha \geq 0$ be a rational number, and let
$u \in {\Bbb R}^d$ be a rational vector. One can show (cf. Section 4.1) that
$$\Big| \bigl((n+\alpha)P + u\bigr) \cap {\Bbb Z}^d \Big|=\sum_{i=0}^d e_i(P, \alpha, u; n) n^i 
\quad \text{for all} \quad n \in {\Bbb N},$$
where $e_i(P, \alpha, u; \cdot): {\Bbb N} \longrightarrow {\Bbb Q}$, $i=0,\ldots, d$, satisfy 
$$e_i(P, \alpha, u; n+t)=e_i(P, \alpha, u; n) \quad \text{for all} \quad n \in {\Bbb N},$$
provided $t \in {\Bbb N}$ is a number such that $tP$ is an integer polytope.
As long as $k$ is fixed in advance, for given $\alpha$, $u$, $n$, and a rational simplex
$\Delta \subset {\Bbb R}^d$, one can compute $e_{d-k}(\Delta, \alpha, u; n)$ in polynomial time.  
\subhead(8.2) Computing the generating function \endsubhead
Let $P \subset {\Bbb R}^d$ be a rational polytope. Then, for every $0 \leq i \leq d$, the 
series 
$$\sum_{n=1}^{+\infty} e_i(P; n)t^n$$
converges to a rational function $f_i(P; t)$ for $|t|<1$.

It is not clear whether $f_{d-k}(\Delta; t)$ can be efficiently computed as 
a ``closed form expression'' in any meaningful sense, although it seems that by adjusting 
the methods of Sections 5--7, for any given $t$ such that $|t|<1$ one can compute the 
value of $f_{d-k}(\Delta; t)$ in polynomial time (again, $k$ is assumed to be fixed in advance).
 \subhead (8.3) Extensions to other classes of polytopes \endsubhead
 If $k$ is fixed in advance, the coefficient $e_{d-k}(P; n)$ can be computed 
 in polynomial time, if the rational polytope $P \subset {\Bbb R}^d$ is given by the 
 list of its $d+c$ vertices or the list of its $d+c$ inequalities, where $c$ is a constant fixed 
 in advance.
 \subhead (8.4) Possible applications to integer programming \endsubhead
 If $P \subset {\Bbb R}^d$ is a rational polytope given by the list of its defining inequalities,
 the problem of testing whether $P \cap {\Bbb Z}^d=\emptyset$ is a typical problem of 
 integer programming, see \cite{G+93} and \cite{Sc86}. Moreover, by a general 
 construction of ``aggregation'' (see Section 16.6 of \cite{Sc86}) the problem can be reduced
 in polynomial time for that for $P=\Delta$. It would be interesting to find out whether 
 efficient computation of $e_{d-k}(\Delta; n)$ can have any practical applications 
 to testing whether $\Delta \cap {\Bbb Z}^d =\emptyset$.
 
 \head References \endhead 
 
 \widestnumber\key{MMMM}

\ref\key{B93a}
\by W. Banaszczyk
\paper  New bounds in some transference theorems in the geometry of numbers
\jour  Math. Ann.
\vol 296 
\yr 1993
\pages 625--635
\endref 

\ref\key{B93b}
\by A. Barvinok
\paper Computing the volume, counting integral points, and exponential sums
\jour  Discrete Comput. Geom. 
\vol 10 
\yr 1993
\pages 123--141
\endref

\ref\key{Ba92}
\by A. Barvinok
\paper Computing the Ehrhart polynomial of a convex lattice polytope
\paperinfo preprint TRITA-MAT-1992-0036
\publ Royal Institute of Technology
\publaddr Stockholm
\yr 1992
\endref

\ref\key{B94a}
\by A. Barvinok 
\paper A polynomial time algorithm for counting integral points in polyhedra when the dimension is 
fixed
\jour Math. Oper. Res.
\vol 19 
\pages  769--779
\yr 1994
\endref

\ref\key{B94b}Ê
\by A. Barvinok
\paper Computing the Ehrhart polynomial of a convex lattice polytope
\jour  Discrete Comput. Geom. 
\vol 12
\yr 1994 
\pages 35--48
\endref

\ref\key{BP99}
\by A. Barvinok and J. Pommersheim
\paper An algorithmic theory of lattice points in polyhedra
\inbook  New perspectives in algebraic combinatorics (Berkeley, CA, 1996--97)
\pages 91--147
\bookinfo Math. Sci. Res. Inst. Publ.
\vol 38 
\publ Cambridge Univ. Press
\publaddr  Cambridge
\yr  1999
\endref

\ref\key{BS05}
\by M. Beck and F. Sottile
\paper Irrational proofs for three theorems of Stanley
\paperinfo preprint arXiv math.CO/0501359
\yr 2005
\endref

\ref\key{Be61}
\by R. Bellman
\book A Brief Introduction to Theta Functions 
\bookinfo Athena Series: Selected Topics in Mathematics 
\publ Holt, Rinehart and Winston
\publaddr  New York 
\yr 1961
\endref

\ref\key{BW03}
\by A. Barvinok and K. Woods 
\paper Short rational generating functions for lattice point problems 
\jour J. Amer. Math. Soc. 
\vol 16 
\yr 2003
\pages  957--979
\endref 

\ref\key{D+04}
\by J. De Loera, R. Hemmecke, M. K\"oppe, and R. Weismantel
\paper Integer polynomial optimization in fixed dimension
\paperinfo preprint arXiv math.OC/0410111
\yr 2004
\endref

\ref\key{DR94}
\by R. Diaz and S. Robins
\paper Solid angles, lattice points, and the Fourier decomposition 
of polytopes
\paperinfo manuscript
\yr 1994
\endref

\ref\key{DR97}
\by R. Diaz and S. Robins 
\paper The Ehrhart polynomial of a lattice polytope 
\jour Ann. of Math. (2)
\vol 145 
\yr 1997
\pages  503--518; Erratum:  146 (1997), no. 1, 237
\endref

\ref\key{GK94}
\by P. Gritzmann and V. Klee
\paper On the complexity of some basic problems in computational convexity. II. Volume and mixed volumes
\inbook Polytopes: Abstract, Convex and Computational (Scarborough, ON, 1993)
\pages 373--466
\bookinfo NATO Adv. Sci. Inst. Ser. C. Math. Phys. Sci.
\vol 440
\publ  Kluwer Acad. Publ.
\publaddr Dordrecht
\yr 1994
\endref

\ref\key{G+93}
\by M. Gr\"otschel, L. Lov\'asz, and A. Schrijver
\book Geometric Algorithms and Combinatorial Optimization. Second edition
\bookinfo  Algorithms and Combinatorics
\vol 2 
\publ Springer-Verlag
\publaddr  Berlin
\yr  1993
\endref

\ref\key{La02}
\by P.D.Lax
\book Functional Analysis
\bookinfo Pure and Applied Mathematics
\publ Wiley-Interscience 
\publaddr New York
\yr 2002
\endref

\ref\key{Ma02}
\by J. Matou\v sek
\book Lectures on Discrete Geometry
\bookinfo Graduate Texts in Mathematics
\vol 212
\publ Springer-Verlag
\publaddr New York
\yr 2002
\endref

\ref\key{Mc78}
\by P. McMullen
\paper Lattice invariant valuations on rational polytopes
\jour  Arch. Math. (Basel)
\vol 31 
\yr 1978/79
\pages 509--516
\endref

\ref\key{Mc93}
\by P. McMullen 
\paper Valuations and dissections 
\inbook Handbook of Convex Geometry
\vol B 
\pages 933--988
\publ North-Holland
\publaddr Amsterdam
\yr 1993
\endref 

\ref\key{Mo93}
\by R. Morelli 
\paper Pick's theorem and the Todd class of a toric variety 
\jour Adv. Math. 
\vol 100 
\yr 1993
\pages 183--231
\endref

\ref\key{MS83}
\by P. McMullen and R. Schneider
\paper Valuations on convex bodies 
\inbook Convexity and its Applications
\pages 170--247 
\publ Birkh\"auser
\publaddr Basel
\yr 1983
\endref

\ref\key{PT04}
\by J. Pommersheim and H. Thomas
\paper Cycles representing the Todd class of a toric variety 
\jour J. Amer. Math. Soc. 
\vol 17 
\yr 2004
\pages 983--994
\endref

\ref\key{Sc86}
\by A. Schrijver
\book Theory of Linear and Integer Programming 
\bookinfo Wiley- \break Interscience Series in Discrete Mathematics
\publ John Wiley $\&$ Sons, Ltd.
\publaddr Chichester
\yr 1986
\endref

\ref\key{Sc93}
\by R. Schneider
\book Convex Bodies: the Brunn-Minkowski Theory
\bookinfo Encyclopedia of Mathematics and its Applications
\vol  44 
\publ Cambridge University Press
\publaddr Cambridge
\yr 1993
\endref

 \ref\key{St97}
 \by R.P. Stanley
 \book Enumerative Combinatorics. Vol. 1
 \bookinfo Corrected reprint of the 1986 original. Cambridge Studies in Advanced Mathematics, 49
 \publ Cambridge University Press
 \publaddr Cambridge
 \yr 1997
 \endref
  \enddocument